\documentclass[12pt]{article}

\usepackage{amsmath,amssymb,dsfont}
\usepackage{graphicx,color}
\usepackage{enumitem}

\newtheorem{theorem}{Theorem}[section]

\def\R{\mathbb R}

\def\SS{{\mathcal S}^N}
\def\xb{{\bar x}}
\def\yb{{\bar y}}

\begin{document}

\title{Local Gradient Estimates for Second-Order Nonlinear Elliptic and Parabolic Equations by the Weak Bernstein's Method }

\author{G.Barles\thanks{Institut Denis Poisson (UMR CNRS 7013)
Université de Tours, Université d'Orléans, CNRS. Parc de Grandmont
 37200 Tours,    France. Email: guy.barles@idpoisson.fr    \newline
    \indent This work was partially supported by the project ANR MFG (ANR-16-CE40-0015-01) funded by the French National Research Agency } }

\maketitle
 \noindent {\bf Key-words}: Second-order elliptic and parabolic equations, gradient bounds, weak Bernstein's method, viscosity solutions.
\\
{\bf MSC}:
35D10  
35D40,   
35J15  
35K10  

\begin{abstract}
{\footnotesize In the theory of second-order, nonlinear elliptic and parabolic equations, obtaining local or global gradient bounds is often a key step for proving the existence of solutions but it may be even more useful in many applications, for example to singular perturbations problems. The classical Bernstein's method is a well-known tool to obtain these bounds but, in most cases, it has the defect of providing only a priori estimates. The  ``weak Bernstein's method'', based on viscosity solutions' theory, is an alternative way to prove the global Lipschitz regularity of solutions together with some estimates but it is not so easy to perform in the case of local bounds. The aim of this paper is to provide an extension of the ``weak Bernstein's method'' which allows to prove local gradient bounds with reasonnable technicalities.}
\end{abstract}

\maketitle

The classical Bernstein's method is a well-known tool for obtaining gradient estimates for solutions of second-order, elliptic and parabolic equations (cf. Caffarelli and Cabr\'e \cite{CClivre} Gilbarg and Trudinger\cite{GT} (Chap. 15)  and Lions\cite{LB}). The underlying idea is very simple: if $\Omega$ is a domain in $\R^N$ and $u : \Omega \to \R$ is a smooth solution of 
$$ -\Delta u = 0 \quad\hbox{in  }\Omega\; ,$$
where $\Delta$ denotes the Laplacian in $\R^N$, then $w:=|Du|^2$ satisfies
$$ -\Delta w \leq 0 \quad\hbox{in  }\Omega\; .$$
The gradient bounded is deduced from this property by using the Maximum Principle if one knows that $Du$ is bounded on $\partial \Omega$ and this bound on the boundary is usually the consequence of the existence of barriers functions.

Of course this strategy, consisting in showing that $w:=|Du|^2$ is a {\em subsolution} of an elliptic equation and then using the
Maximum Principle, can be applied to far more general equations but it has a clear defect: in order to justify the above
computations, the solution has to be $C^3$ and, since it is rare that the solution has such a regularity, the classical Bernstein's method provides, in general, only {\em a priori estimates}; then one has to find a suitable approximation of the equation, with smooth enough solutions, to actually obtain the gradient bound.

In 1990, this difficulty was partially overcomed by the weak Bernstein's method whose idea is even simpler: if one looks at the maximum of the function
$$(x,y) \mapsto u(x)-u(y)-L|x-y|  \quad\hbox{in  }\overline \Omega \times \overline \Omega\; ,$$
and if one can prove that it is achieved only for $x=y$ for $L$ large enough, then $|Du|\leq L$. Surprisingly, as it is explained in the introduction of \cite{B-wb}, the computations and structure conditions which are needed to obtain this bound are the same (or almost the same with tiny differences) as for the classical Bernstein's method. Of course, the main advantage of the weak Bernstein's method is that it does not require $u$ to be smooth since there is no differentiation of $u$ and it can even be used in the framework of viscosity solutions.

Problem solved? Not completely because the weak Bernstein's method is not of an easy use if one looks for local bounds instead of global bounds. In fact, in order to get such local gradient bounds, the only possible way seems to multiply the solution by a cut-off function and to look for a gradient bound for this new function. Unfortunately, this new function satisfies a rather complicated equation where the derivatives of the cut-off function appear at different places and the computations become rather technical. The classical Bernstein's method also faces similar difficulties but, at least in some cases, succeeds in providing these local bounds in a not too complicated way.

The aim of this article is to describe a slight improvement of the weak Bernstein's method which allows to obtain local gradient bounds in a simpler way, ``simpler'' meaning that the technicalities are as reduced as possible, although some are unavoidable. This improvement is based on an idea of P.~Cardaliaguet~\cite{C1} which dramatically simplifies a matrix analysis which is keystone in \cite{B-wb} but also allows this extension to local bounds.

To present our result, we consider second-order, possibly degenerate, elliptic equations which we write in the general form
\begin{equation}\label{GFNLE}
F(x, u, D u, D^2u) = 0  \quad\hbox{in  }\Omega\; ,
\end{equation}
where $\Omega$ is a domain of $\R^N$ and $F :\Omega \times \R \times \R^N \times \SS \to \R $ is a locally Lipschitz continuous function, $\SS$ denotes the space of $N \times N$ symmetric matrices, the solution $u$ is a real-valued function defined on $\Omega$, $Du, D^2u$ denote respectively its gradient and Hessian matrix. We assume that $F$ satisfies the (degenerate) ellipticity condition : for any $(x,r,p)\in\Omega \times \R \times \R^N$ and for any $X,Y\in\SS$, 
$$
F(x,r,p,X) \leq F(x,r,p,Y)\quad\hbox{if }X\geq Y.
$$

Our results consist in providing several general ``structure conditions'' on $F$ under which one has a local gradient bound
depending or not on the local oscillation of $u$ and the uniform ellipticity of the equation. We also consider the parabolic case
for which we give a structure condition on the equation allowing to prove a local gradient bound, depending on the local oscillation of
$u$, where ``local'' means both in space and time.

In the stationary framework, we focus in particular on the following example
\begin{equation}\label{PartEqn}
-\Delta u + |Du|^m = f(x)  \quad\hbox{in  }\Omega\; ,
\end{equation}
where $m>1$ and $f \in W^{1,\infty}_{loc}(\Omega)$, which is a particular case for which the classical Bernstein's method
provides local bound (independent of the oscillation of $u$)  in a rather easy way, while it is not the case for the weak Bernstein's method.

We conclude this introduction by two remarks: the first one concerns the ``structure conditions'' on $F$ on which our results are based.
In \cite{B-wb}, it is pointed out that, in general, the equation we consider does not satisfy these structure conditions and we have to make a 
change of unknown function $v=\psi(u)$, choosing $\psi$ in order that the new equation for $v$ satisfies them. Obviously, 
the same remark is true here and we provide an example where such a change allows to obtain the desired gradient bound. But, contrarily to 
\cite{B-wb}, we are not going to study the effect of such changes in a more systematic way.

The second remark concerns the method we are going to present: the results we obtain are based on several choices we made at several
places and, in particular, in the estimates of the terms we have to handle. Clearly, many variants are possible and we have just tried to
 convince the reader that, 
actually, the technicalities are really ``reasonnable'' as we pretend it in the abstract.\\

\noindent{\bf Acknowledgement:} the author would like to thank the anonymous referees whose remarks led to significant improvements of the readability of this article.

\section{Some preliminary results}\label{prelim}

In this section, we are going to construct the functions we use in the proof of our main result. To do so, we introduce $\mathcal{K}$ which is the class of continuous functions $\chi:[0,+\infty)\to [0,+\infty)$ such that $\chi(t)=0$ if $t\leq 1$, $\chi$ is increasing on $[1,+\infty[$, $\chi(t)\leq {\tilde K}(\chi)t^\beta$ for $t\geq 1$, for some $0<\beta < 1/2$ and some constant ${\tilde K}(\chi)>0$, and
$$ \int_1^{+\infty}\frac{dt}{t\chi(t)}<+\infty .$$

The first ingredient we use below is a smooth function $\varphi : [0,1[ \to \R$ such that $\varphi(0)=0$, $\varphi'(0)=1 \leq \varphi'(t)$ for any $t\in [0,1[$ with $\varphi(t) \to +\infty$ as $t\to 1^-$ and which solves the ode $\varphi''(t)= K_1\varphi'(t) \chi(\varphi'(t))$ for some constant $K_1>0$. In fact the existence of such function is classical using that
$$ \int_1^{\varphi'(t)} \frac{ds}{s\chi(s)} = K_1 t\, ,$$
and by choosing $K_1=\int_1^{+\infty} \frac{ds}{s\chi(s)}$ we already see that $\varphi' (t) \to +\infty$ as $t\to 1^-$. Moreover 
$$ \int_{\varphi'(t)}^{+\infty} \frac{ds}{s\chi(s)} = K_1 (1-t) \; ,$$
and therefore, for $t$ close enough to $1$
$$ K_1 (1-t) \geq [{\tilde K}(\chi)]^{-1}\int_{\varphi'(t)}^{+\infty} \frac{ds}{s^{1+\beta}}=  [{\tilde K}(\chi)\beta]^{-1}\varphi'(t)^{-\beta}\; .$$
This means that
$$\varphi'(t) \geq \left(\frac{K_1 (1-t)}{[{\tilde K}(\chi)\beta]^{-1}}\right)^{-1/\beta} \; ,$$
and therefore $\varphi'(t)$ is not integrable at $1$ since $1/\beta>2$. Hence we have $\varphi(t) \to +\infty$ as $t\to 1^-$.

On the other hand, given $x_0 \in \R^N$ and $R>0$, we use below a smooth function $C: B(x_0,3R/4) \to \R$ is a smooth function such that $C(z)= 1$ on $B(x_0,R/4)$, $C(z) \geq 1$ in $ B(x_0,3R/4)$ and $C(z)\to +\infty$ when $z\to \partial B(x_0,3R/4)$ and with
$$ \frac{|D^2C(x)|}{C(x)} , \frac{|DC(x)|^2}{[C(x)]^2} \leq K_2(R) [\chi(C(x))]^2\; ,$$ where $\chi$ is a function in the class $\mathcal{K}$. If $C_1$ is a function which satisfies the above properties for $x_0=0$ and $R=1$, we see that we can choose $C$ as
$$ C(x)=C_1\left(\frac{x-x_0} R\right)\; ,$$
and therefore $K_2(R)$ behaves like $R^{-2}K_2(1)$.

To build $C_1$, we first solve
$$ \psi'' (t) = K_3 \psi (t)[\chi(\psi (t))]^2, \; \psi(0)=1,\;  \psi'(0)=0\; ,$$
for some constant $K_3$ to be chosen later on.
Multiplying the equation by $2 \psi'(t)$, we obtain that
$$ \psi'(t) = F(\psi(t))\; ,$$
where
$$ [F(\tau)]^2= 2K_3 \int_1^\tau s[\chi(s)]^2 ds\; .$$
Again we look for a function $\psi$ such that $\psi(t) \to +\infty$ as $t\to 1^{-}$ and to do so, the following condition should hold
$$ \int_1^{+\infty} \frac{d\tau}{F(\tau)} < +\infty\; .$$
But, since $\chi$ is increasing,
$$ [F(\tau)]^2 \geq 2K_3 \int_{\tau/2}^\tau s[\chi(s)]^2 ds\geq \; 2K_3 [\tau/2 \chi(\tau/2)]^2 ,$$
and since $\tau \mapsto \chi(\tau/2)$ is in $\mathcal{K}$, we have the result for $F$, and then for $\psi$ by choosing appropriately the constant $K_3$.

Moreover 
$$ [F(\tau)]^2 \leq  2K_3 (\tau-1) \tau[\chi(\tau)]^2 \leq 2K_3 [\tau \chi(\tau)]^2 \; ,$$
and therefore
$$ \psi'(t) \leq (2K_3)^{1/2} \psi(t) \chi(\psi(t))\; .$$
Finally, we can extend $\psi$ by setting $\psi(t)=1$ for $t\leq 0$ and the equations satisfied by $\psi$ show that we define in that way a $C^2$-function on $(-\infty,1)$.

With such a $\psi$, the construction of $C_1$ is easy, we may choose
$$ C_1(x):= \psi \bigl(4(|x| - 1/2)\bigr)\quad \hbox{for  }x\in B(0,3/4) ,$$
and define $C$ from $C_1$ as above. We notice that, because of the properties of $\psi$, $\dfrac{|DC(x)|}{[C(x)]^2}$ remains bounded on
$B(x_0,3R/4)$ and is a $O(R^{-1})$, a property that we will use later on.
 
\section{The Main Result}

In the statement of our main result below, for the sake of clarity, we are going to drop the arguments of the partial derivatives of $F$ and
to simply denote by $F_s$ the quantity $\dfrac{\partial F}{\partial s} (x,r,p,M)$ for $s=x,r,p,M$. Actually these arguments are $(x,r,p,M)$ everywhere.

Our result is the following

\begin{theorem} \label{main}Assume that $F$ is a locally Lipschitz function in $\Omega \times \R \times \R^N \times \SS \to \R$ which satisfies : $F(x,r,p,M)$ is Lipschitz continuous in $M$ and
$$
F_M(x,r,p,M) \leq 0 \;\hbox{and}\; F_r(x,r,p,M) \geq 0\quad\hbox{a.e. in  }\Omega \times \R \times \R^N \times \SS\; ,$$
and let $u\in C(\Omega)$ be a solution of (\ref{GFNLE}).\\

(i) {\bf (Uniformly elliptic equation with coercive gradient dependence: estimates which are independant of the oscillation of $u$)} Assume that there exist a function $\chi \in \mathcal{K}$ and $0<\eta\leq1$ such that, for any $K>0$, there exists $L= L(F,K)$ large enough such that
$$ -(1+\eta)|F_x| |p| (1+K\chi(\eta |p|)) - K |F_p|  |p|^2 \left(1+K\chi(\eta |p| )\right) \chi(\eta |p| ) - \dfrac1{1+\eta}F_M\cdot M^2 $$
$$ \geq \eta + K \bigl( |p| \left(1+K\chi(\eta |p| )\right) \chi(\eta |p| )\bigr)^2 \; \hbox{a.e.},
$$
in the set $$\{(x,r,p,M);\ |F(x,r,p,M))| \leq K \eta |p|[1 + K\chi(\eta|p|)]+\eta\; ,\; |p|\geq L\}\; .$$
If $\overline{B(x_0,R)} \subset \Omega$ then $u$ is Lipschitz continuous in $B(x_0,R/2)$ and $|Du| \leq {\bar L} $ in $B(x_0,R/2)$ where $\bar L$ depends only on $F$ and $R$.\\

(ii) {\bf (Uniformly elliptic equation with coercive gradient dependence: estimates depending the oscillation of $u$)} Assume that there exist a function $\chi \in \mathcal{K}$ and $0<\eta\leq1$ small enough such that, for any $K>0$, there exists $L= L(F,K)$ large enough such that
$$-(1+\eta) |F_x||p| - K|F_p| |p|^2\chi(\eta |p|) - \frac{1}{1+\eta} F_M\cdot M^2 \geq \eta 
+  K |p|^2\chi(\eta |p|)^{2}\; \hbox{a.e.},$$
in the set $\{(x,r,p,M);\ |F(x,r,p,M))| \leq K |p|+\eta \; ,\; |p|\geq L\}$. If $\overline{B(x_0,R)} \subset \Omega$ then $u$ is Lipschitz continuous in 
$B(x_0,R/2)$ and $|Du| \leq {\bar L} $ in $B(x_0,R/2)$ where $\bar L$ depends on $F$, $R$ and $osc_R (u)$, the oscillation of $u$ on $
\overline{B(x_0,R)}$.\\

(iii) {\bf (Non-uniformly elliptic equation : estimates depending the oscillation of $u$)} Assume that there exist a function $\chi \in \mathcal{K}$
and $0< \eta \leq 1$ small enough such that, for any $K>0$, there exists $L=L(F,K)$ large enough such that
$$ -(1+\eta) |F_x| |p| +(1-\eta)^2 F_r|p|^2 - K|F_p| |p|^2\chi(\eta |p|)- \frac{1}{1+\eta} F_M\cdot M^2$$
$$ \geq \eta +  K |p|^2\chi(\eta |p|)^{2}\; \hbox{a.e.},$$
in the set $\{(x,r,p,M);\ |F(x,r,p,M))| \leq K |p|+\eta \; ,\; |p|\geq L\}$. If $\overline{B(x_0,R)} \subset \Omega$ then $u$ is Lipschitz continuous in $B(x_0,R/2)$ and $|Du| \leq {\bar L} $ in $B(x_0,R/2)$ where $\bar L$ depends on $F$, $R$ and $osc_R (u)$.
\end{theorem}

{ As an application we consider Equation (\ref{PartEqn}): in order to have a gradient estimate which is independant of the oscillation of $u$, i.e. Result (i) in Theorem~\ref{main},
 the idea is to choose $\chi(t)=(t-1)^\beta$ for $t\geq 1$ with $0<\beta <1/2$ and $\gamma:=1+2\beta < m$. The most important point is that, for large $|p|$, the constraint on $F$ reads
 $$|F(x,r,p,M))| \leq K\eta |p|(1+ K(\eta |p|)^{\beta})+\eta$$ and therefore $|F(x,r,p,M))|$ behaves as $K^2(\eta |p|)^{1+\beta}$ if $|p|$ is
 large enough. Since $1+\beta<m$, this implies that, for such $(x,r,p,M)$, 
$$ {\rm Tr}(M)\geq \frac12 |p|^m - ||f||_{L^{\infty}(B(x_0,R)}\; .$$
But, by Cauchy-Schwarz inequality
$$ {\rm Tr}(M)\leq C(N)[{\rm Tr}(M^2)]^{1/2}\; .$$
Therefore the term $-F_M\cdot M^2$ behaves like $|p|^{2m}$. For the other terms, we have, for large $|p|$
\begin{enumerate}
\item the term $|F_x|  |p| (1+K\chi(\eta |p|)) $ behaves like $|p|^{1+\beta}=|p|^{\gamma-\beta}$;
\item the term $|F_p|   |p|^2 \left(1+K\chi(\eta |p| )\right) \chi(\eta |p| )$ behaves like $|p|^{m+1+2\beta}=|p|^{m+\gamma}$;
\item the term $K ||F_M||_\infty \bigl( |p| \left(1+K\chi(\eta |p| )\right) \chi(\eta |p| )\bigr) ^2$ behaves like $|p|^{2(1+2\beta)}=|p|^{2\gamma}$.
\end{enumerate}
Since $\gamma < m$, the term $-F_M\cdot M^2$ clearly dominates all the other terms as $|p|$ tends to $+\infty$; therefore we have the gradient bound since the assumption holds for any $0<\eta \leq 1$. Moreover the classical case ($m=1$) can be also treated under the assumptions of Result~(ii).

In this example, it is also clear that we can replace the term $|Du|^m$ by a term $H(Du)$ where $H$ satisfies: there exists $\chi \in  \mathcal{K}$ such that
$$\frac{|p|\chi(|p|)}{H(p)}\to 0 \quad\hbox{as  } |p|\to +\infty\; , $$
and 
$$ \frac{|H_p|(|p|\chi(|p|))^2}{[H(p)]^2}\to 0 \quad\hbox{as  } |p|\to +\infty\; .$$}
\ \\

In the case of non-uniformly elliptic equation, the gradient bound comes necessarely from the $F_r|p|^2$-term. We consider the equation
\begin{equation}\label{PartEqnNUN}
-{\rm Tr}(A(x)D^2 u) + |Du|^m = f(x)  \quad\hbox{in  }\Omega\; ,
\end{equation}
where $m>1$ and $f $ is locally bounded and Lipschitz continuous; concerning $A$, we use the classical assumption: $A(x)=\sigma(x)\cdot \sigma^T(x)$ for some bounded, Lipschitz continuous function $\sigma$, where $\sigma^T(x)$ denotes the transpose matrix of $\sigma(x)$.

In order to obtain a local gradient bound for $u$, a change of variable is necessary: assuming (without loss of generality) that $u\geq 1$ at least in the ball $\overline{B(x_0,R)}$, we can use the change $u=\exp(v)$. The equation satisfied by $v$ is
$$
-{\rm Tr}(A(x)D^2 v) +A(x)Dv\cdot Dv+ \exp((m-1)v)|Dv|^m = \exp(-v)f(x)  \quad\hbox{in  }\Omega\; ,
$$
And the aim is now to apply Theorem~\ref{main}-(iii) to get the gradient bound for $v$ (hence for $u$).

The computation of the different terms gives
$$ F_r(x,r,p,M)= (m-1)\exp((m-1)r)|p|^m + \exp(-r)f(x)\; ,$$
$$ F_x(x,r,p,M)= -{\rm Tr}(A_x(x)M)+A_x(x)p\cdot p-\exp(-r)f_x (x)$$
$$F_p(x,r,p,M)= 2A(x)p +\exp((m-1)r)|p|^{m-2}p\; ,$$
$$ - F_M(x,r,p,M)M^2= {\rm Tr}(A(x)M^2)\; .$$

We first use Cauchy-Schwarz inequality and the assumption on $A$ to deduce that, for any $\eta>0$
$$ |{\rm Tr}(A_x(x) M)|  |p|\leq \frac{1}{1+\eta} {\rm Tr}(A(x)M^2)+ O((|\sigma_x||p|)^2)\; ;$$
This control of the first term in $F_x(x,v,p,M)$ is the only use of the term $- F_M(x,v,p,M)M^2$ .

Therefore the $F_r(x,r,p,M)|p|^2$-term which behaves like $|p|^{m+2}$ if $m>1$, has to control the terms 
$$ (A_x(x)\cdot p)(p\cdot p)=O(|p|^3)\; ,\; -\exp(-v)f_x (x) |p|=O(|p|)\; ,\; 2A(x)p\cdot p =O(|p|^2)\; .$$

We have now to consider the $F_p$-term and the term $K \bigl(|p|\chi (\eta |p|)\bigr)^{2}$ in the right-hand side. Notice that, for the time being, we have not chosen $\chi$ nor $\eta$.

The $F_p$-term behaves as $|p|^{\max(1,m-1)}$ and therefore $|F_p|   |p|^2\chi(\eta |p|)$ behaves as $|p|^{\max(3,m+1)}\chi(\eta |p|)$.
On the other hand, $K \bigl(|p|\chi (\eta |p|)\bigr)^{2}$ behaves as $|p|^2[\chi(\eta |p|)]^2$. If we choose any $\chi \in \mathcal{K}$, because of the growth of such $\chi$ at infinity, these two terms are controlled by the $F_r|p|^2$-one. Therefore Theorem~\ref{main} (iii) applies.

It is worth pointing out that, in this last example, we do not use the fact that the assumption has to hold only
in the set $\{(x,r,p,M);\ |F(x,r,p,M))| \leq K |p|+\eta \; ,\; |p|\geq \bar L\}$, a fact which is going to be (almost) the general case in the parabolic setting.
 {
\section{Proof of Theorem~\ref{main}}

We start by proving (i) : the aim is to prove that, for any $x\in B(x_0,R/4)$, $D^+u(x)$ is bounded with an explicit bound. This will provide the desired gradient bound. We recall that $$ D^+u(x)=\{p\in\R^n:\  u(x+h)\leq u(x)+p\cdot h+o(|h|) \ \text{ as }  h\to 0\}.$$ 

To do so, we consider on 
$$\Gamma_L :=\{(x,y) \in B(x_0,3R/4) \times B(x_0,R) : LC(x)(|x-y| +\alpha)<1\}$$ 
the following function
$$ \Phi(x,y)= u(x)-u(y) - \varphi\left (LC(x)(|x-y| +\alpha)\right )\; ,$$
where
\begin{itemize}
\item $L\geq \max(1,4/R)$ is a constant which is our future gradient bound (and therefore which has to be choosen large enough),
\item the functions $\varphi$ and $C$ are built in Section~\ref{prelim},
\item $\alpha >0$ is a small constant devoted to tend to $0$.
\end{itemize}

We remark that the above function achieves its maximum in the open set $\Gamma_L$: indeed, if $(x,y) \in \Gamma_L$, we have $ LC(x)\alpha<1$ and therefore $x\in \overline{B(x_0,R')}$ for some $R'<3R/4$. Moreover $LC(x)|x-y|<1$ implies $|x-y|<L^{-1}$ and, since $L> 4/R$, this implies $y\in \overline{B(x_0,R'+R/4)}$ and $R'+R/4<R$. Therefore, clearly $\Phi(x,y) \to - \infty$ if $(x,y)\to \partial \Gamma_L$.

Next we argue by contradiction: if, for some $L$, this maximum is achieved for any $\alpha$ at $(\xb_\alpha,\yb_\alpha)$ with $\xb_\alpha=\yb_\alpha$, then $\Phi(\xb_\alpha,\xb_\alpha)=- \varphi(LC(\xb_\alpha)\alpha)$ and therefore necessarely $\xb_\alpha \in B(x_0,R/4)$ by the maximality property and the form of $C$. Moreover, for any $x,y$
$$u(x)-u(y) - \varphi(LC(x)(|x-y| +\alpha))\leq -\varphi(L\alpha)\; ,$$
and if this is true, for a fixed $L$, this implies that, for any $x,y$
$$u(x)-u(y) - \varphi(LC(x)|x-y| )\leq 0\; .$$
Choosing $x\in B(x_0,R/4)$, we have
$$ u(y)-u(x) \geq - \varphi(L|x-y| )\; ,$$
and this inequality implies that any element in $D^+u(x)$ has a norm which is less than $L$, which we wanted to prove.

Notice that, by using slightly more complicated arguments, the same conclusion is true if, for some $L$, we have $\xb_\alpha-\yb_\alpha \to 0$ when $\alpha \to 0$.

Therefore, we may assume without loss of generality that, for any fixed $L$, the maximum points $(\xb_\alpha,\yb_\alpha)$ of $\Phi$, satisfies not only $\xb_\alpha\neq \yb_\alpha$ for $\alpha$ small enough but $\xb_\alpha-\yb_\alpha$ is bounded away from $0$ when $\alpha \to 0$. We are going to prove that this is a contradiction for $L$ large enough.

For the sake of simplicity of notations, we omit the indice $\alpha$ in all the quantities which depends on $\alpha$ (actually they also depend on $L$). In particular, we denote by $(x,y)$ a maximum point of $\Phi$ and we set $t=LC(x)(|x-y| +\alpha)$ and
$$ p= \varphi'(t)LC(x) \frac{(x-y)}{|x-y|} \; ,\; q=  \varphi'(t)LDC(x) (|x-y|+\alpha)\; .$$
By a classical result of the User's guide (cf. Crandall, Ishii and Lions \cite{users}), there exist matrices $X,Y \in \SS$ such that $(p+q,X)\in \overline{D^{2,+}}u(x)$, $(p,Y)\in \overline{D^{2,-}}u(y)$, for which the following viscosity inequalities hold
$$ F(x,u(x), p+q,X) \leq 0\; ,\; F(y,u(y), p,Y) \geq 0\; .$$
Moreover the matrices $X,Y$ satisfy, for any $\varepsilon >0$
$$ \left(-\frac1\varepsilon + ||A||\right) I_{2N} \leq \left(\begin{array}{cc}X & 0 \\0 & -Y\end{array}\right)\leq A + \varepsilon A^2$$
and where, if $\psi(x,y)= \varphi(LC(x)(|x-y| +\alpha))$, $A=D^2\psi(x,y)$ and $||A||=\max\{|\lambda|:\ \hbox{$\lambda $ is an eigenvalue of $A$}\}$.

Since $\varepsilon>0$ is arbitrary and since we are going to use only the second above inequality, we may choose a sufficiently
small $\varepsilon$ in order that the term $\varepsilon A^2$ becomes negligible. Using this remark, we argue
below assuming that $\varepsilon=0$ in order to simplify the exposure.

With this convention, the matrices $X,Y$ satisfy, for any $r,s \in \R^N$
\begin{equation}\label{fu}
Xr\cdot r - Ys \cdot s \leq \gamma_1|r-s|^2+2\gamma_2 |r-s||r|+\gamma_3|r|^2\; ,
\end{equation}
where
$$ \gamma_1= \frac{\varphi'(t)LC(x)}{|x-y|}+ \varphi''(t)(LC(x))^2\; ,$$
$$ \gamma_2= \varphi'(t)L|DC(x)|+ \varphi''(t)L^2|DC(x)|C(x)(|x-y|+\alpha)\; ,$$
$$ \gamma_3= \varphi'(t)\frac{|D^2C(x)|}{C(x)}t+ \varphi''(t)\frac{|DC(x)|^2}{[C(x)]^2}t^2\; ,$$
By easy manipulations, it is easy to see that
$$ \gamma_2 \leq \gamma_1\frac{|DC(x)|}{C(x)}(|x-y|+\alpha) + o_\alpha(1)\leq \gamma_1 K_2^{1/2}\chi(C(x))(|x-y|+\alpha)+ o_\alpha(1)\; ,$$
$$ \gamma_3 \leq \gamma_1 K_2[\chi(C(x))]^2(|x-y|+\alpha)^2+ o_\alpha(1)\; ,$$
where the $o_\alpha(1)$ comes from terms of the form $\alpha/|x-y|$. Again, for the sake of clarity, we are going to drop these terms which play no role at the end.

By Cauchy-Schwarz inequality, we deduce that, using $\eta$ appearing in the assumption,
\begin{equation}\label{si}
Xr\cdot r - Ys \cdot s \leq (1+\eta)\gamma_1|r-s|^2+ B(R,\eta)\gamma_1 [\chi(C(x))]^2(|x-y|+\alpha)^2|r|^2\; ,
\end{equation}
where $B(R,\eta)=(1+\eta^{-1})K_2$ depends on $R$ through $K_2$ and therefore is a $O(R^{-2})$ if $\eta$ is fixed.

Coming back to $p$ and $q$, we also have
$$ |q| =  |p|\frac{|DC(x)|}{C(x)} (|x-y|+\alpha) \leq  |p|\frac{|DC(x)|}{L [C(x)]^2} \leq O((RL)^{-1})|p|\; ,$$
since $LC(x)(|x-y|+\alpha)\leq 1$, $C\geq 1$ everywhere and since $\dfrac{|DC(x)|}{ [C(x)]^2}$ is a $O(R^{-1})$.
In order to have simpler formulas, we denote below by $\varpi_1$ any quantity which is a $O((RL)^{-1})$.

Now we arrive at the key point of the proof: by \eqref{fu}, choosing $r=0$, we have $-Y\leq \gamma_1I_N$ where $I_N$ is the identity matrix in $\R^N$. Therefore the matrix $\displaystyle I_N+[(1+\eta)\gamma_1]^{-1}Y$ is invertible and rewriting \eqref{si} as
$$ Xr\cdot r  \leq  Ys \cdot s + (1+\eta)\gamma_1|r-s|^2+ B(R,\eta)\gamma_1 [\chi(C(x))]^2(|x-y|+\alpha)^2|r|^2 \; ,$$
we can take the infimum in $s$ in the right-hand side and we end up with 
$$X \leq Y(I_N+\frac{1}{(1+\eta)\gamma_1}Y)^{-1}+B(R,\eta)\gamma_1 [\chi(C(x))]^2(|x-y|+\alpha)^2I_N\; .$$

Setting  $\tilde Y:= Y(I_N+\frac{1}{(1+\eta)\gamma_1}Y)^{-1}$, this implies that we have $(p+q,\tilde Y+3\gamma_1 [\chi(C(x))]^2(|x-y|+\alpha)^2I_N)\in \overline{D^{2,+}}u(x)$, $(p,Y)\in \overline{D^{2,-}}u(y)$ and then, using the Lipschitz continuity of $F$ in $M$, we have
the viscosity inequalities
$$ F(x,u(x), p+q,\tilde Y) \leq ||F_M||_\infty B(R,\eta)\gamma_1 [\chi(C(x))]^2(|x-y|+\alpha)^2\; ,$$
$$  F(y,u(y), p,Y) \geq 0\;.$$

Next we introduce the function
$$
g(\tau):= F(X(\tau), U(\tau), P(\tau), Z(\tau))-\tau||F_M||_\infty B(R,\eta)\gamma_1 [\chi(C(x))]^2(|x-y|+\alpha)^2 \; ,$$
where 
$$ X(\tau) = \tau x+(1-\tau)y\; ,\; U(\tau)= \tau u(x)+(1-\tau)u(y)\; ,\; P(\tau)=p+\tau q\; ,$$
$$ Z(\tau) = Y(I_N+\frac{\tau}{(1+\eta)\gamma_1}Y)^{-1}\; .$$

From now on, in order to simplify the exposure, we are going to argue as if $F$ were $C^1$: the case when $F$ is just locally Lipschitz continuous follows from tedious but standard approximation arguments.

The above viscosity inequalities read $g(0)\geq 0$ and $g(1)\leq 0$ : if we can show that the $C^1$-function $g$ satisfies $g'(\tau)>0$ if $g(\tau)=0$, we would have a contradiction. Therefore we compute
\begin{eqnarray*}
g'(\tau) &=& F_x\cdot(x-y)+F_r (u(x)-u(y)) + F_p\cdot q + F_M\cdot Z'(\tau)\\
&&-||F_M||_\infty B(R,\eta)\gamma_1 [\chi(C(x))]^2(|x-y|+\alpha)^2\; ,
\end{eqnarray*}
and using that $F_r \geq 0$ and $Z'(\tau)=-((1+\eta)\gamma_1)^{-1}[Z(\tau)]^2$, we are lead to
\begin{eqnarray*}
g'(\tau) &\geq &  (\gamma_1)^{-1}\biggl\{-|F_x|   \gamma_1 |x-y|- \gamma_1 |F_p|   |q| - \dfrac1{1+\eta}F_M\cdot [Z(\tau)]^2\\
&&  -B(R,\eta) ||F_M||_\infty (\gamma_1)^2 [\chi(C(x))]^2(|x-y|+\alpha)^2\biggr\}.
\end{eqnarray*}

Before estimating the different terms inside the brackets, we point out that, contrarily to \cite{B-wb} where $Z(\tau)$ was given by
$\tau X + (1-\tau)Y$ and where we had to prove an inequality between $X-Y$ and $-[Z(\tau)]^2$, here this inequality comes for free
because of the form of $Z(\tau)$: this is the key idea of Cardaliaguet \cite{C1}.

Now we estimate the terms $\gamma_1 |x-y|$, $\gamma_1 |q|$ and $\gamma_1 \chi(C(x))(|x-y|+\alpha)$ in terms of $|P(\tau)|$ in
order to be able to use the assumptions on $F$.

Using that $LC(x)(|x-y|+\alpha) \leq 1$ and the properties of $\varphi$, we have
\begin{eqnarray*}
\gamma_1|x-y| & \leq & \varphi'(t)LC(x)+ \varphi''(t)(LC(x))^2|x-y|\\
& \leq & \varphi'(t)LC(x) + K_1\varphi'(t)\chi(\varphi'(t)) LC(x)\\
& \leq & |P(\tau)| (1+\varpi_1)\left(1 + K_1\chi(\varphi'(t)) \right)\\
& \leq & |P(\tau)| (1+\varpi_1)\left(1+K_1 \chi(L^{-1}|P(\tau)| (1+\varpi_1))\right)\; .
\end{eqnarray*}
Indeed, recalling the estimate on $|q|$, $ \varphi'(t)LC(x)=|p|=|P(\tau)|(1+\varpi_1\tau)$ and, on an other hand, since $C\geq 1$, we have
$$\chi(\varphi'(t))\leq \chi(L^{-1}|p|)\leq \chi(L^{-1}|P(\tau)| (1+\varpi_1)).$$

From now on, we are going to assume that $L$ is chosen large enough in order to have $L^{-1} (1+\varpi_1)\leq \eta$ and, since $R$ is fixed, $|\varpi_1| \leq \eta$. Notice that these constraints on $L$ depend only on $R$ and $\eta$, hence on $R$ and $F$.

Using this choice, the above estimate of $\chi(\varphi'(t))$ -- and we can argue in the same way for $\chi(C(x))$-- takes the simple form
\begin{equation}\label{kest}
\chi(\varphi'(t)), \chi(C(x)) \leq \chi(\eta |P(\tau)|).
\end{equation}
This leads to the
simpler estimate
$$ \gamma_1|x-y|  \leq |P(\tau)| (1+\eta)(1+K_1\chi(\eta |P(\tau)|))\; .$$
In the same way, since we can take $\alpha$ as small as we want and $|x-y|$ is bounded away from $0$, one has
$$ \gamma_1 (|x-y| +\alpha) \leq |P(\tau)| (1+\eta)\left(1+K_1\chi(\eta |P(\tau)|)\right) +o_\alpha(1)\; .$$
This allows to estimate the $F_p$-term, namely
\begin{align*}
\gamma_1 |q| & \leq \gamma_1 |p| \frac{|DC(x)|}{C(x)}(|x-y| +\alpha)\\
&\leq |P(\tau)|^2 (1+\eta)^2 \left(1+K_1\chi(\eta |P(\tau)| )\right)K_2^{1/2} \chi(C(x))+o_\alpha(1),\\
& \leq K_2^{1/2} |P(\tau)|^2 (1+\eta)^2 \left(1+K_1\chi(\eta |P(\tau)| )\right) \chi(\eta |P(\tau)| )+o_\alpha(1)\; .
\end{align*}
Finally, by the same estimates
$$\gamma_1 \chi(C(x)) (|x-y|+\alpha) \leq |P(\tau)| (1+\eta)\left(1+K_1\chi(\eta |P(\tau)| )\right)\chi(\eta |P(\tau)| ) +o_\alpha(1)\; .$$

We end up with
\begin{eqnarray*}
g'(\tau) &\geq &  (\gamma_1)^{-1}\biggl\{-|F_x|   |P(\tau)| (1+\eta)(1+K_1\chi(\eta |P(\tau)|) ) \\
&& \phantom{(\gamma_1)^{-1}\biggl\{}  - |F_p| K_2^{1/2} |P(\tau)|^2 (1+\eta)^2 \left(1+K_1\chi(\eta |P(\tau)| )\right) \chi(\eta |P(\tau)| ) \\
&& \phantom{(\gamma_1)^{-1}\biggl\{} - \dfrac1{1+\eta}F_M\cdot [Z(\tau)]^2\\
&&   -B(R,\eta) ||F_M||_\infty \bigl( |P(\tau)| (1+\eta)\left(1+K_1\chi(\eta |P(\tau)| )\right) \chi(\eta |P(\tau)| )\bigr) ^2\biggr\} \\
&& + o_\alpha(1).
\end{eqnarray*}
On the other hand, in order to take into account the constraint $g(\tau)= 0$, we have to estimate $\gamma_1 [\chi(C(x))]^2 (|x-y|+\alpha)^2$. 
Since $|x-y|$ is bounded away from $0$ and $LC(x)(|x-y|+\alpha)\leq 1$, we have
\begin{eqnarray*}
\gamma_1(|x-y|+\alpha)^2 & \leq & \varphi'(t)+ \varphi''(t) +o_\alpha (1)\\
& \leq & \varphi'(t)[1 + K_1\chi(\varphi'(t))]+ o_\alpha (1)\\
& \leq & (1+\eta) \frac{|P(\tau)|}{LC(x)}[1 + K_1\chi(\varphi'(t))]+o_\alpha (1)\\
& \leq & (1+\eta) \frac{|P(\tau)|}{LC(x)}[1 + K_1\chi(\eta|P(\tau)| )]+o_\alpha (1).
\end{eqnarray*}
But $\dfrac{[\chi(C(x))]^2}{C(x)} \leq \tilde K(\chi)$ and therefore
$$ \gamma_1 [\chi(C(x))]^2 (|x-y|+\alpha)^2 \leq \eta(1+\eta)\tilde K(\chi) |P(\tau)|[1 + K_1\chi(\eta |P(\tau)| )]+o_\alpha (1).$$
This implies
$$ |F(X(\tau), U(\tau), P(\tau), Z(\tau))|\leq  \eta(1+\eta)\tilde K(\chi) |P(\tau)|[1 + K_1\chi(\eta|P(\tau)|)]+o_\alpha (1)\; ,
$$
while
$$ |P(\tau)| \geq (1-\eta)L\; .$$

The conclusion follows by applying the assumption on $F$ for $L$ large enough and $\alpha$ small enough in order that the $o_\alpha(1)$-terms are controlled by the $\eta$-terms. Taking $L$ large enough depending on $\eta$ and $R$, we have a contradiction and the proof of (i) is complete.

Now we turn to the proof of (ii) where we choose $\varphi(t)=t$ and 
$$\Gamma'_L :=\{(x,y) \in B(x_0,3R/4) \times B(x_0,R) : LC(x)(|x-y| +\alpha)\leq osc_R (u)\}\; .$$
The proof follows the same arguments, except that the fact that $\varphi''(t)\equiv 0$ allows different estimates on the $\gamma_i$, $i=1,2,3$ because several terms do not exist anymore. We denote by $\varpi_2$ any quantity of the form $O(osc_R (u)(RL)^{-1})$ and we choose $L$ 
large enough in order to have $|\varpi_2|\leq \eta$ for any of these terms and $L^{-1} \leq \eta/(1+\eta)$. We notice that, here, the constraints on $L$ depend not only on $R$ and $\eta$ but also on $osc_R (u)$.

We have $p= LC(x) \dfrac{(x-y)}{|x-y|}$ and therefore
$$ |q|=  L.|DC(x)| (|x-y|+\alpha)= |p|   \frac{|DC(x)|}{C^2}\frac{LC(x)(|x-y|+\alpha)}{L}=\varpi_2|p|\leq \eta |p|\; ,$$
since $\dfrac{|DC(x)|}{C^2}\leq O(R^{-1})$. Using this inequality and taking into account our choice of $L$, it is easy to check that (\ref{kest}) still holds. 

Moreover we have
$$ \gamma_1= \frac{LC(x)}{|x-y|}\; ,\;  \gamma_2= L.|DC(x)|\; ,
\; \gamma_3= L |D^2C(x)|(|x-y|+\alpha) \; .$$
And we still have the same estimates on $\gamma_1, \gamma_2,\gamma_3$
$$ \gamma_2 = \gamma_1\frac{|DC(x)|}{C(x)}|x-y| \leq \gamma_1 \chi(C(x))|x-y|\; ,$$
$$ \gamma_3 \leq \gamma_1 [\chi(C(x))]^2(|x-y|+\alpha)^2\; .$$

The proof is then done in the same way as in the first case with the computation of $ g'(\tau)$ and then with the estimates of the different terms
\begin{eqnarray*}
g'(\tau) &\geq &  (\gamma_1)^{-1}\left\{-|F_x|   \gamma_1 |x-y|-\gamma_1 |F_p|   |q| - \frac{1}{1+\eta}F_M\cdot [Z(\tau)]^2\right.\\
&& \left . -B(R,\eta) ||F_M||_\infty (\gamma_1)^2 [\chi(C(x))]^2(|x-y|+\alpha)^2\right\}\; .
\end{eqnarray*}
But here 
$$\gamma_1 |x-y|=|p| \leq |P(\tau)| (1+\eta)\; ,$$
and in the same way,
\begin{align*}
\gamma_1 |q| & = \frac{LC}{|x-y|} L |DC(x)|(|x-y| +\alpha)\\
&\leq |p|^2 \dfrac{|DC(x)|}{C(x)}(1+o_\alpha(1))\\
& \leq K_2^{1/2} (1+\eta)^2 |P(\tau)|^2  \chi(\eta |P(\tau)| )+o_\alpha(1)\; ,
\end{align*}
and
$$\gamma_1 \chi(C(x)) (|x-y|+\alpha) \leq (1+\eta)|P(\tau)|\chi(\eta|P(\tau)|)+o_\alpha (1)\; .$$

We end up with
\begin{eqnarray*}
g'(\tau) &\geq &  (\gamma_1)^{-1}\biggl\{-|F_x|   (1+\eta)|P(\tau)| - K_2^{1/2} (1+\eta)^2|F_p|   |P(\tau)|^2  \chi(\eta |P(\tau)| ) \\
&& -\frac{1}{1+\eta} F_M\cdot [Z(\tau)]^2  - B(R,\eta) ||F_M||_\infty (1+\eta)^2 |P(\tau)|^2[\chi(\eta|P(\tau)|)]^2 \\
&& +o_\alpha (1) \biggr\}\; .
\end{eqnarray*}
On the other hand, for the constraint $g(\tau)= 0$, we have
\begin{align*}
\gamma_1[\chi(C(x))]^2(|x-y|+\alpha)^2 &= |p|\frac{[\chi(C(x))]^2}{C(x)}\dfrac{LC(|x-y|+\alpha)^2}{|x-y|}\\
 &\leq (1+\eta)[\tilde K(\chi)]^2 |P(\tau)|(1+\varpi_2)(1+o_\alpha(1))\\
&\leq  (1+\eta)^2[\tilde K(\chi)]^2|P(\tau)|+o_\alpha(1)\; ,
\end{align*}
and 
$$ |P(\tau)|\geq LC(x) (1-\eta)\geq L(1-\eta)\; .$$
Hence
\begin{equation}\label{F-prop}
 |F(X(\tau), U(\tau), P(\tau), Z(\tau))| \leq B(R,\eta)||F_M||_\infty (1+\eta)^2 |P(\tau)| + o_\alpha(1)\; .
\end{equation}
The conclusion follows as in the first case by applying the assumption on $F$ for $L$ large enough and $\alpha$ small enough for which we have a contradiction.

For the proof of (iii), we keep the same test-function and the same set $\Gamma'_L$ but since we are not expecting the gradient bound to come from the same term in $g'(\tau)$, we are going to change the strategy in our computation of $g'(\tau)$ by keeping the $F_r$-term. Using that
$F_r \geq 0$ and 
$$ u(x)-u(y)\geq LC(x)(|x-y|+\alpha)= \frac{|p |^2}{\gamma_1}(1+o_\alpha(1))\; ,$$
we obtain
\begin{eqnarray*}
g'(\tau)& = & F_x\cdot(x-y)+F_r (u(x)-u(y)) + F_p\cdot q + F_M\cdot Z'(\tau)\\
&&-||F_M||_\infty B(R,\eta)\gamma_1 [\chi(C(x))]^2(|x-y|+\alpha)^2\; ,\\
 &\geq&  (\gamma_1)^{-1}\left\{F_x\cdot p+ F_r |p|^2 - \gamma_1 |F_p|   |q| - \frac{1}{1+\eta}F_M\cdot [Z(\tau)]^2\right.\\
&& \left . -B(R,\eta) ||F_M||_\infty (\gamma_1)^2 [\chi(C(x))]^2(|x-y|+\alpha)^2+o_\alpha(1)\right\}\; .
\end{eqnarray*}
This computation is close to the one given in \cite{B-wb} if there is no localization term ($C\equiv1$).

Since $|P(\tau)| (1-\eta)\leq |p| \leq |P(\tau)| (1+\eta)$ and using anagolous estimates as above, we are lead to
\begin{eqnarray*}
g'(\tau) &\geq &  (\gamma_1)^{-1}\biggl\{-(1+\eta) |F_x|    |P(\tau)| +(1-\eta)^2 F_r |P(\tau)|^2\\
&& - K_2^{1/2} (1+\eta)^2|F_p|   |P(\tau)|^2  \chi(\eta |P(\tau)| )- \frac{1}{1+\eta}F_M\cdot [Z(\tau)]^2 \\ 
\\ 
&&  -B(R,\eta) ||F_M||_\infty (1+\eta)^2 [\tilde K(\chi)]^2 |P(\tau)|^2[\chi(\eta|P(\tau)|)]^2 \biggr\}+o_\alpha (1) \; .
\end{eqnarray*}
On the other hand, the constraint $g(\tau)= 0$ still implies (\ref{F-prop}) and we also conclude by choosing $L$ large enough and $\alpha$ small enough.}

\section{The parabolic case}

In this section, we consider evolution equations under the general form
\begin{equation}\label{GFNLP}
u_t + F(x, t, u, D u, D^2u) = 0  \quad\hbox{in  }\Omega \times (0,T)\; ,
\end{equation}
and the aim is to provide a local gradient bound where ``local'' means both local in space and time. As a consequence, we will have to provide a localization also in time and a second main difference is that we will not be able to use that the equation holds since the $u_t$-term has no property in general and therefore the assumptions on $F$ have to hold for any $x, t, r, p, M$ and not only those for which $F(x, t, r, p, M)$ is close to $0$.

\begin{theorem} \label{mainP}{\bf (Estimates for non-uniformly parabolic equations : estimates depending the oscillation of $u$)}\\
Assume that $F$ is a locally Lipschitz function in $\Omega \times (0,T) \times \R \times \R^N \times \SS$ which satisfies : $F(x,t,r,p,M)$ is Lipschitz continuous in $M$ and 
$$
F_M(x,t,r,p,M) \leq 0 \quad\hbox{a.e. in  }\Omega \times (0,T)\times \R \times \R^N \times \SS \; ,$$
and let $u\in C(\Omega\times (0,T))$ be a solution of (\ref{GFNLP}). Assume that there exists a function $\chi \in \mathcal{K}$, $0<\eta\leq 1$
 such that, for any $K>0$, there exists $L=L(\eta,K)$ large enough such that, for $|p|\geq  L$,
we have $F_r(x,t,r,p,M)\geq 0$ and
$$ -(1+\eta) |F_x|   |p| (1+\chi(\eta |p|)) - K |F_p|  |p|^2 \left(1+\chi(\eta |p| )\right) \chi(\eta |p| )- \frac{1}{1+\eta} F_M\cdot M^2$$
$$ \geq \eta +  K |p|^2\biggl( \chi((1+\eta) |p|)+ \chi(\eta |p|)^{2}\biggr)\; \hbox{a.e.},$$
If $\overline{B(x_0,R)} \subset \Omega$ and $\delta>0$, then $u$ is Lipschitz continuous in $x$ in $B(x_0,R/2)\times [\delta, T-\delta]$ and $|Du| \leq {\bar L} $ in $B(x_0,R/2)\times [\delta, T-\delta]$ where $\bar L$ depends on $F$, $R$, $\delta$ and the oscillation of $u$ in
$B(x_0,R)\times (\delta/2,T-\delta]$.
\end{theorem}

It is worth pointing out that the assumptions of Theorem~\ref{mainP} are rather close to the one of Theorem~\ref{main} (iii) and the same computations provide a gradient bound for the evolution equation
\begin{equation}\label{PartEqnNUN-P}
u_t-{\rm Tr}(A(x)D^2 u) + |Du|^m = f(x)  \quad\hbox{in  }\Omega\times (0,T)\; ,
\end{equation}
if $m>1$.

\noindent{\bf Proof of Theorem~\ref{mainP} :} We argue as in the proof of Theorem~\ref{main} (iii), except that here $L=L(t)$ with $L(t) \to +\infty$ as $t \to (\delta/2)^+$. We still choose $\varphi(t)=t$ and we denote by $\Gamma'_L$, the subset of points $(x,y,t) \in B(x_0,3R/4) \times B(x_0,R)\times (\delta/2,T-\delta]$ such that
$$ L(t)C(x)(|x-y| +\alpha)\leq osc_{R,\delta} (u)\},$$
where $osc_{R,\delta} (u)$ denotes the oscillation of $u$ in
$B(x_0,R)\times (\delta/2,T-\delta]$.

We consider maximum points $(x,y,t) \in \Gamma'_L$ of the function
$$(x,y,t)\mapsto u(x,t)-u(y,t) - L(t)C(x)(|x-y| +\alpha)\; ,$$
and, if $x\neq y$, we are lead to the viscosity inequalities
$$ a+ F(x,t, u(x,t), p+q,X) \leq 0\; ,\; b+F(y,t,u(y,t), p,Y) \geq 0\; ,$$
where $(a,p+q,X)\in D^{2,+}u(x,t)$, $(p,Y)\in D^{2,-}u(y,t)$ and $$a-b \geq L'(t)C(x)(|x-y|+\alpha).$$

As in the proof of Theorem~\ref{main}, the second inequality holds for $\tilde Y$ as well and subtracting these inequalities, we have
$$ L'(t)C(x)(|x-y|+\alpha)+ F(x,u(x), p+q,X) -F(y,u(y), p,\tilde Y)\leq 0\; .$$
Then, with the notations of the proof of Theorem~\ref{main}, we introduce
$$
 g(\tau) :=  F(X(\tau), U(\tau), P(\tau), Z(\tau))-\tau ||F_M||_\infty B(R,\eta)\gamma_1 [\chi(C(x))]^2(|x-y|+\alpha)^2]$$
$$ +\tau L'(t)C(x)(|x-y|+\alpha)  \; .$$
Here we have no information on the signs of $g(0)$ and $g(1)$, we only know that $g(1)-g(0)\leq 0$; therefore, in order to have the contradiction, we have to show that $g'(\tau) >0$ for any $0\leq\tau\leq1$ if we choose a function $L(\cdot)$ such that $L(t)$ is large enough for any $t\in (\delta/2,T-\delta]$.

The computation of $g'(\tau)$ and the estimates are done as above; we have just to estimate the new term $L'(t)C(x)(|x-y|+\alpha)$ which is
multiplied by $\gamma_1$ when we put it inside the bracket. We have
$$ \gamma_1 L'(t)C(x)(|x-y|+\alpha) = L(t) L'(t)[C(x)]^2 (1+o_\alpha(1)) \; ,$$
and if we choose $L$ as the solution of the ode
$$ L'(t)=-k_T L(t)\chi(L(t))\; , L(T-\delta) = L_T \;\hbox{(large enough)}\; .$$
By choosing properly $k_T>0$, we have $L((\delta/2)^+)=+\infty$ (notice that $k_T$ decreases when $L_T$ increases). Since $L(t)\leq |p|\leq (1+\eta) |P(\tau)|$, we have
$$L(t) L'(t)[C(x)]^2 \geq -k_T |P(\tau)|^2 \chi((1+\eta)|P(\tau)|)\; .$$

Using this estimate, the conclusion follows as above by applying the assumption on $F$ for $K$ large enough and $\alpha$ small enough for which we have a contradiction by taking $L_T $ large enough.


\end{document}